\documentclass[11pt]{amsart}

\usepackage{epsfig}
\newcommand{\lcr}{\raisebox{-5pt}{\mbox{}\hspace{1pt}
                  \epsfig{file=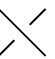}\hspace{1pt}\mbox{}}}
\newcommand{\ift}{\raisebox{-5pt}{\mbox{}\hspace{1pt}
                  \epsfig{file=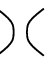}\hspace{1pt}\mbox{}}}
\newcommand{\zer}{\raisebox{-5pt}{\mbox{}\hspace{1pt}
                  \epsfig{file=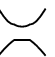}\hspace{1pt}\mbox{}}}

\theoremstyle{definition}

\theoremstyle{remark}

\numberwithin{equation}{section}

\begin{document}

\setlength{\textheight}{7.7truein}  

\title[Noncommutative A-ideal of torus knots]{\bf THE NONCOMMUTATIVE  A-IDEAL
OF A  $(2,2p+1)$-TORUS  KNOT DETERMINES ITS JONES POLYNOMIAL}

\maketitle

\vspace*{0.37truein}
\centerline{\footnotesize R\u{A}ZVAN GELCA and JEREMY SAIN
}
\baselineskip=12pt
\centerline{\footnotesize\it Department of Mathematics and Statistics}
\baselineskip=10pt
\centerline{\footnotesize\it Texas Tech University }
\centerline{\footnotesize\it Lubbock, TX 79409}


\begin{abstract}
{The noncommutative A-ideal of a knot is a generalization of
the A-polynomial, defined using Kauffman bracket skein modules. 
In this paper we show that any knot that has the same noncommutative
A-ideal as the $(2,2p+1)$-torus knot has the same colored Jones polynomials.
This is a consequence of the orthogonality relation, which yields a
recursive relation for computing all  colored Jones polynomials of 
the knot.}
\end{abstract}



\section{Introduction}	

The noncommutative A-ideal was defined in \cite{FGL} using Kauffman
bracket skein modules. Skein modules were introduced by Turaev (for
surfaces) \cite{T1} and Przytycki (for arbitrary 3-manifolds) \cite{P} 
in an attempt to generalize the polynomial invariants of knots and links
in the 3-sphere to invariants of knots and links in arbitrary manifolds.
The simplest skein module is the Kauffman bracket skein module,
defined  using   the Kauffman bracket skein relation
\cite{K}. The module depends on a parameter $t$, 
the variable of the Kauffman bracket.
A  theorem of  Bullock, Przytycki and Sikora
\cite{B1}, \cite{PS} shows  that
when $t=-1$, the Kauffman bracket skein module of a 3-manifold 
has a ring structure, which ring is isomorphic with the affine
character ring of
 $SL(2,{\mathbb C})$-representations of the
fundamental group of the manifold. 

This result allows the study of  $SL(2,{\mathbb C})$-invariants
of knots  from the point of view of quantum invariants. 
An example for which this is done  is  the A-polynomial.
The A-polynomial of a knot was introduced by Cooper, Culler, Gillet, Long,
and Shalen in \cite{CCGLS}  using the  character variety of 
$SL(2,{\mathbb C})$-representations of the fundamental group 
of the knot complement. The
 construction was later  generalized via Kauffman bracket
skein modules in \cite{FGL}, giving rise to a new knot invariant, the 
noncommutative A-ideal. The A-ideal is a finitely generated ideal
of polynomials in two noncommuting variables. It
  depends on a parameter $t$, and
when $t=-1$, the radical of the one dimensional part of the ideal is 
generated by the A-polynomial of the knot. 

In \cite{FGL} was discovered a  property relating
the noncommutative A-ideal to the Jones polynomial.
 Called the orthogonality relation,
this property states that every nonzero element in the A-ideal yields a
matrix annihilating the vector with entries  the colored Jones
polynomials of the knot. As an application of the orthogonality relation,
it was shown in \cite{G2} that if the A-ideal contains a polynomial of
second degree in the variable  corresponding to the 
longitude of the knot, whose coefficients satisfy a certain technical property,
then the A-ideal of the knot determines all of its colored Jones polynomials.
In this paper we prove that this is the case with the $(2,2p+1)$-torus knots, 
$p$ an integer.
The main result of the paper is 

{\bf Theorem~1.1} {\it Any knot having  the same noncommutative A-ideal as the $(2,2p+1)$-torus 
knot has the same colored Jones polynomials.}

To avoid the hassle of keeping track of both the sign of $p$ and 
of the types of crossings in  diagrams, we do the proof 
for the case where $p$ is {\em positive}. The case where $p$ is negative
is obtained by replacing in all formulas $p$ by $|p|$ and 
$t$ by $t^{-1}$.  

\section{Preliminary Facts}

A {\em framed link} in an orientable manifold $M$ is a 
disjoint union of annuli. If the manifold is
the cylinder over the torus, framed links will be identified with 
curves, using the convention that the annulus is parallel to the framing.
In figures   annuli will  be represented by curves, the framing
will be  considered  parallel to the plane of the drawing (the so called 
blackboard framing). 
Consider the free ${\mathbb C}[t]$-module ${\mathbb C}[t]{\mathcal L}$
with basis $\mathcal{L}$
the set of all isotopy classes of links in $M$, including the empty link.
The quotient of this module by the smallest submodule containing
all expressions
of the form $\displaystyle{\lcr-t\zer-t^{-1}\ift}$
and 
$\bigcirc+t^2+t^{-2}$, where the links in each expression are
identical except in a ball in which they look like depicted,
is called the Kauffman bracket skein module of $M$, and is
denoted by $K_t(M)$. If $t$ is a complex number instead of
the variable of a polynomial, then with the same definition the
Kauffman bracket skein module becomes a complex vector space.

The Kauffman bracket skein module of a cylinder over a surface has
an algebra structure. The multiplication is induced by the operation
of gluing one cylinder atop another. Similarly, 
the operation of gluing the cylinder over $\partial M$ to $M$ induces
a $K_t(\partial M\times I)$-left module structure  on $K_t(M)$.
We denote by $*$ the multiplication in $K_t(\partial M\times I)$.
For a link $\gamma $ in a skein module we will denote by
$\gamma ^n$ the link consisting of $n$ parallel copies of 
$\gamma$, and extend the notation to polynomials.

 Two families of polynomials
are of interest to us. The first family 
consists of 
 the classical Chebyshev polynomials, $T_0=2$, $T_1=x$, $T_{n+1}=xT_n-
T_{n-1}$. Related to them are the polynomials $S_n$ subject to the same
recurrence relation but with $S_0=1$, $S_1=x$. 
Extend both polynomials recursively  to
all indices $n\in{\mathbb Z}$. Note that $T_{-n}=T_n$, while
$S_{-n}=-S_{n-2}$. For a knot $K$, $S_n(K)$ is called the coloring
of $K$ by the $n$th Jones-Wenzl idempotent. If $K$ is a knot
in the 3-sphere, then $S_n(K)$ as a polynomial in $K_t(S^3)={\mathbb C}[t]$ 
is called the $n$th colored Kauffman bracket  of the knot. Under the change of
variable $t\rightarrow it$ this becomes the $n$th colored Jones
polynomial. Although this is a small alteration, and so 
$S_n(K)$ itself could very well be called the colored Jones polynomial,
we prefer the name colored Kauffman bracket, to keep track of the fact
that it comes from the Kauffman bracket. As such, we denote it by
$\kappa _n(K)$.

Let us describe now $K_t({\mathbb T}^2\times I)$,
the Kauffman bracket skein algebra of the cylinder over
a torus. Let $p,q$ be  two integers, with $n$ their common divisor and
$p'=p/n$, $q'=q/n$. We denote by   $(p,q)_T$ be the skein $T_n((p',q'))$ in 
  $K_t({\mathbb T}^2\times I)$, where $(p',q')$ is the simple closed curve
of slope $p'/q'$ on the torus.
 As a module, $K_t({\mathbb T}^2\times I)$ is free with basis $(p,q)_T$, 
$p\geq 0$, $q\in {\mathbb Z}$. As shown  in \cite{FG}, the multiplication
is given by the product-to-sum formula
\begin{eqnarray*}
(p,q)_T*(r,s)_T=t^{|^{pq}_{rs}|}
(p+r, q+s)_T+
t^{-|^{pq}_{rs}|}
(p-r,q-s)_T.
\end{eqnarray*}
This  formula shows that there is an inclusion of 
$K_t(M)M$ into the ring of trigonometric functions on the noncommutative torus.
Recall that this ring, denoted by 
${\mathbb C}_t[l,l^{-1},m,m^{-1}]$, consists of Laurent polynomials in $l$ and
$m$, where the two variables satisfy $lm=t^2ml$. The inclusion map is 
defined by 
\begin{eqnarray*}
(p,q)_T\rightarrow t^{-pq}(l^pm^q+l^{-p}m^{-q}).
\end{eqnarray*} 
A subring  of  ${\mathbb C}_t[l,l^{-1},m,m^{-1}]$ is the quantum plane
${\mathbb C}_t[l,m]$, consisting of the polynomials in $l$ and $m$, again
with $lm=t^2ml$. 

 If $M$ is the complement of a regular neighborhood of a knot $K$, then
the noncommutative A-ideal of $K$, ${\mathcal A}_t(K)$ is defined
as follows.  Denote by $\pi$ the map between skein modules induced
by the inclusion $\partial M \times I\subset M$ and let $I_t(K)$ be the
kernel of $\pi$. $I_t(K)$ is a left ideal called the peripheral ideal
of K. The noncommutative A-ideal of $K$ is defined to 
be the left ideal obtained by extending $I_t(K)$ to 
${\mathbb C}_t[l,l^{-1},m,m^{-1}]$,
then contracting it to ${\mathbb C}_t[l,m]$.

\setcounter{footnote}{0}
\renewcommand{\thefootnote}{\alph{footnote}}

\section{The Kauffman Bracket Skein Module of the Complement
of a $(2,2p+1)$-Torus Knot}

Denote by $M_p$ the complement of
a regular neighborhood of the $(2,2p+1)$-torus knot. 
Bullock proved in  \cite{B2} that $K_t(M_p)$ is a free module with basis 
\begin{eqnarray*}
\{x^ky^n; \quad  0\leq k, 0\leq n\leq p\},
\end{eqnarray*}
where $x$  and $y$  are depicted in Fig. 3.1.  
Here $x^ky^n$ means $k$ parallel copies of $x$ together with $n$ parallel
copies of $y$. It is important to remark that although the
skein module of the knot complement does not have a multiplicative operation,
this formula makes sense. 

\begin{figure}[htbp]
\centering
\leavevmode
\epsfxsize=2.2in
\epsfysize=1.2in
\epsfbox{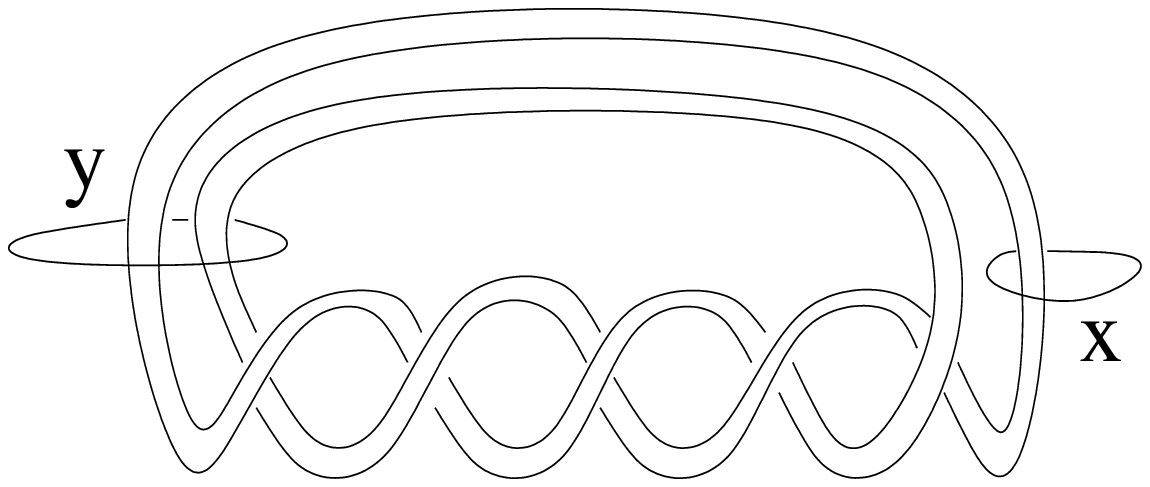}

Figure 3.1.    
\end{figure}

Computations are simpler if we work instead with the basis
\begin{eqnarray*}
\{ S_k(x)S_n(y);\quad 0\leq k, 0\leq n\leq p\}.
\end{eqnarray*}

The following result about the Kauffman bracket skein module of
the complement of the $(2,2p+1)$-torus knot is important in itself. In addition
it  will be used extensively in our computations. 

{\bf Theorem~3.1.} {\it
For  $i=0, 1, \ldots, p+1$, one has
\begin{eqnarray*}
t^{-2i-1}S_{p+i}(y) +t^{2i+1}S_{p-i-1}(y)
 =(-1)^iS_{2i}(x)(tS_{p-1}(y)+t^{-1}S_p(y)).
\end{eqnarray*}}

In order to prove the theorem we 
 first  study the  skeins $A(k,n)$ and $\bar{A}(k,n)$ in
the knot complement  defined 
in Fig. 3.2. It will be seen below that $A(k,n)=\bar{A}(k,n)$ and that
$A(n,k)$ is a polynomial of degree $k+n$ in $y$. Note that
$A(i,0)=\bar{A}(2p+1-i,0)$, which will produce the relation from the 
statement of the theorem. 
\begin{figure}[htbp]
\centering
\leavevmode
\epsfxsize=4.8in
\epsfysize=1.4in
\epsfbox{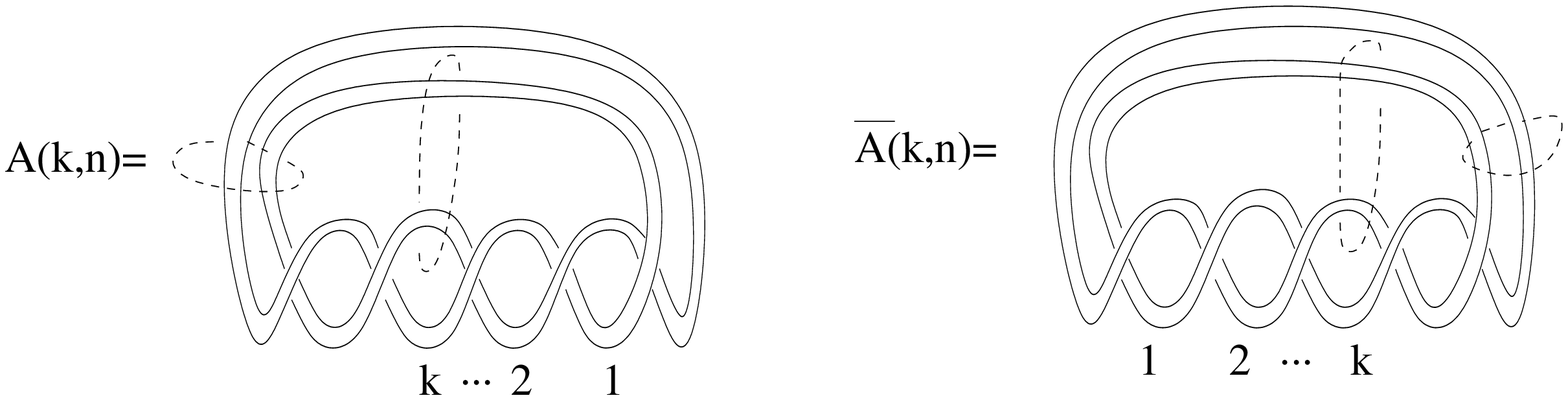}

Figure 3.2.    
\end{figure}

We now prove two  recursive relations.

{\bf Lemma~3.2.} {\it
For all $n\geq 0$ and all $1\leq k\leq 2p$ one has
\begin{eqnarray*}
& & A(k+1,n)=t^{-2}A(k,n+1)-t^{-4}A(k-1,n)+t^2x^2y^n-t^{-2}x^2y^n+x^2y^{n+1}\\
& & \bar{A}(k+1,n)=t^{-2}\bar{A}(k,n+1)-t^{-4}\bar{A}(k-1,n)+t^2x^2y^n-
t^{-2}x^2y^n
+x^2y^{n+1}.
\end{eqnarray*}
}

\begin{proof}
We discuss only the first recursive relation, the second is similar
and is left to the reader. The computation of $A(k+1,n)$ 
begins as described in Fig.3.3.
\begin{figure}[htbp]
\centering
\leavevmode
\epsfxsize=5.2in
\epsfysize=2.5in
\epsfbox{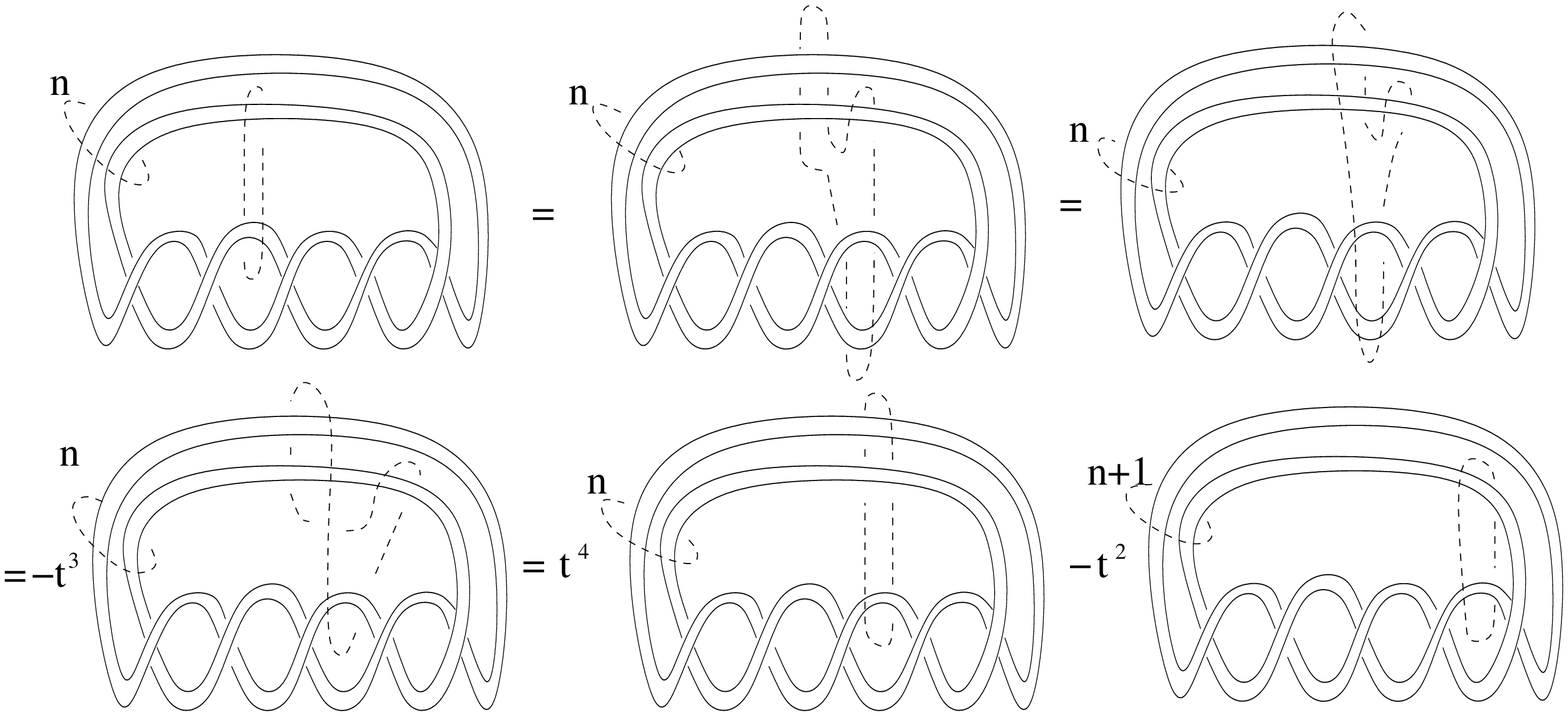}

Figure 3.3.    
\end{figure}

After introducing a ``kink'' and resolving the crossing the
second skein becomes $x^2y^{n+1}+t^{-2}A(k,n+1)$. The first skein is
computed as in Fig. 3.4. 

\begin{figure}[htbp]
\centering
\leavevmode
\epsfxsize=5.2in
\epsfysize=1.3in
\epsfbox{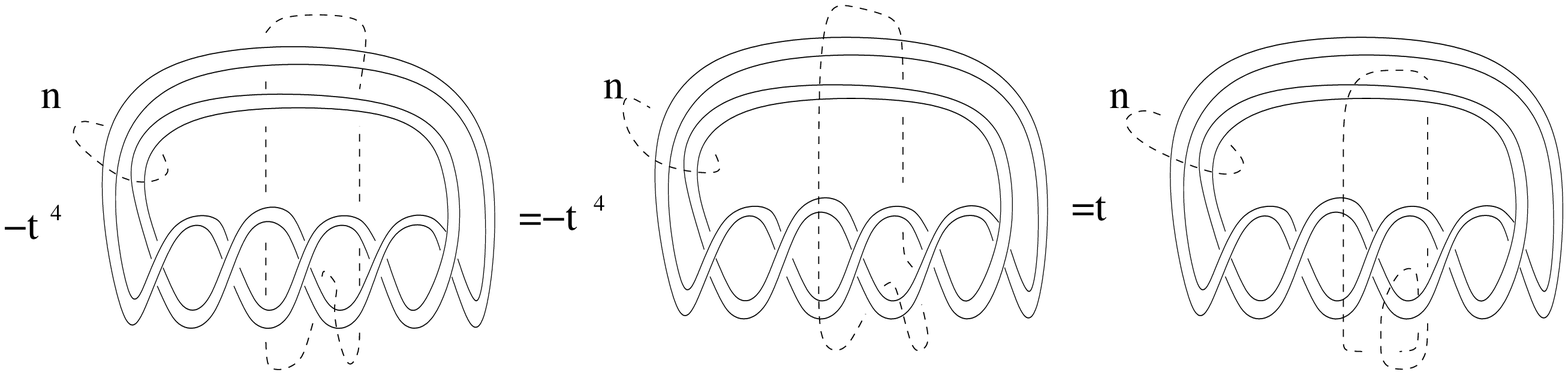}

Figure 3.4.    
\end{figure}

Resolving the crossing we obtain $t^2x^2y^n$ plus the skein from
Fig. 3.5. After introducing one twist and resolving it
we obtain $-t^{-2}x^2y^n-t^{-4}A(k-1,n)$ and the conclusion
of the lemma follows.
\begin{figure}[htbp]
\centering
\leavevmode
\epsfxsize=2in
\epsfysize=1.1in
\epsfbox{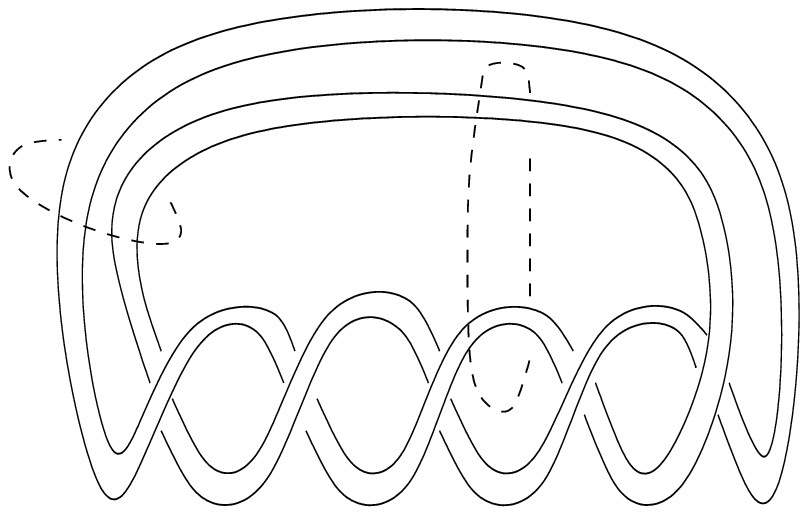}

Figure 3.5.    
\end{figure}
\end{proof}

Now we determine the initial conditions of the two recurrences.

{\bf Lemma~3.3.} {\it
For any $n\geq 0$ one has
\begin{eqnarray*}
& & A(1,n)=\bar{A}(1,n)=y^{n+1}\\
& & A(2,n)=\bar{A}(2,n)=-(t^2+t^{-2})y^n+t^{-2}y^{n+2}+t^2x^2y^n+x^2y^{n+1}.
\end{eqnarray*}
}

\begin{proof}
The case of $A(1,n)$ and $\bar{A}(1,n)$ is easy and is left to the 
reader. 
The computation  of $A(2,n)$ begins 
 as in Fig. 3.6. The computation of the remaining skein is described
in Fig 3.7, and the desired formula is obtained by resolving the two
crossings. A similar computation for $\bar{A}(2,n)$ shows that
$A(2,n)=\bar{A}(2,n)$. In fact one could show directly this equality
by sliding strands. 
 \begin{figure}[htbp]
\centering
\leavevmode
\epsfxsize=5.2in
\epsfysize=1.2in
\epsfbox{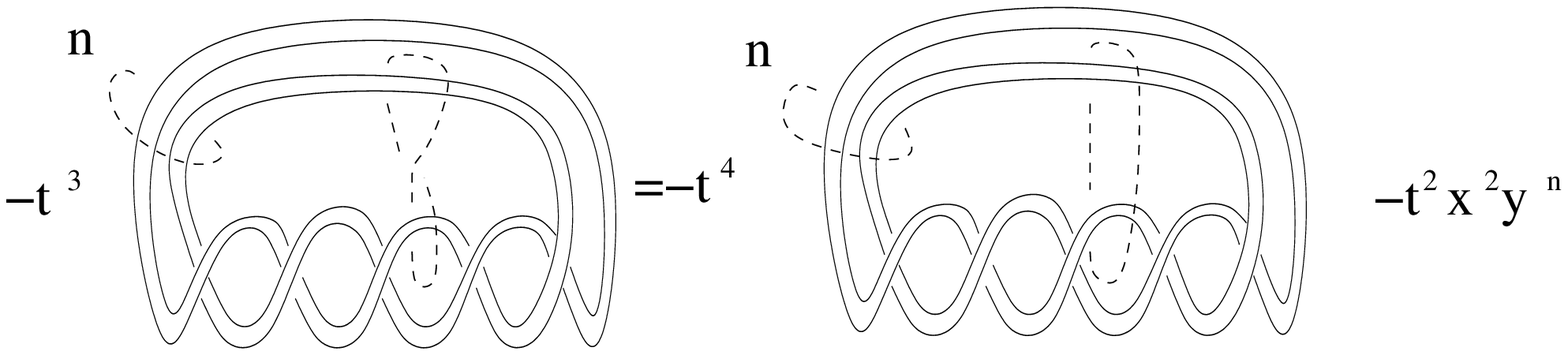}

Figure 3.6. 
 \centering
\leavevmode
\epsfxsize=5.2in
\epsfysize=1.2in
\epsfbox{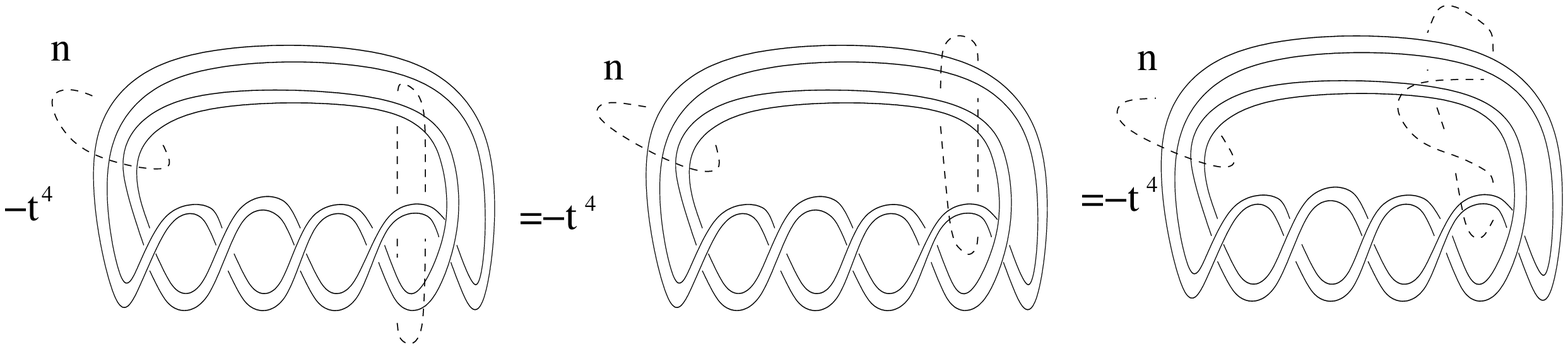}

Figure 3.7.   
\end{figure}
\end{proof}

{\bf Lemma~3.4.} {\it
For all $n\geq 0$ and $1\leq k\leq 2p$ we have
\begin{eqnarray*}
A(k,n)& = & -t^{-2k+6}y^nS_{k-2}(y)-t^{-2k+4}x^2y^nS_{k-1}(d)+
t^{-2k+2}y^nS_k(d)\\
& & -t^2x^2y^n+2t^2x^2\sum_{r=-1}^{k-1}t^{-2r}S_r(y)y^n.
\end{eqnarray*}
}

\begin{proof}
The proof is by induction on $k$. The formula is valid for $k=1$ and
$2$ by Lemma 3.3. Assume that it holds for $k-1$ and $k$ and let us prove
it for $k+1$. Applying  Lemma 3.2 and the induction hypothesis we have
\begin{eqnarray*}
& & A(k+1,n)=-t^{-2k+4}y^{n+1}S_{k-2}(y)-t^{-2k+2}x^2y^{n+1}S_{k-1}(y)\\
& & \quad +
t^{-2k}y^{n+1}S_{k}(y)-x^2y^{n+1}
+2x^2\sum_{r=0}^{k-1}t^{-2r}S_r(y)y^{n+1}\\
& & \quad -(-t^{-2k+4}y^nS_{k-3}(y)
-t^{-2k+2}x^2y^nS_{k-2}(y)+t^{-2k}y^nS_{k-1}(y)\\
& & \quad +2x^2\sum_{r=1}^{k-1}t^{-2r}
S_{r-1}(y)y^n-t^{-2}x^2y^n)
+t^2x^2y^n-t^{-2}x^2y^n+x^2y^{n+1}\\
& & \quad  =-t^{-2k+4}y^nS_{k-1}(y)-t^{-2k+2}x^2y^nS_k(y)+
t^{-2k}y^nS_{k+1}(y)+2x^2y^{n+1}\\
& & \quad +t^2x^2y^n+
\sum_{r=1}^{k-1}2t^{-2r}x^2y^nS_{r+1}(y)
\end{eqnarray*}
\begin{eqnarray*}
& & \quad =-t^{-2k+4}y^nS_{k-1}(y)-t^{-2k+2}x^2y^nS_k(y)+
t^{-2k}y^nS_{k+1}(y)+t^2x^2y^n\\
& & \quad +2x^2y^{n+1}+2t^2x^2\sum_{r=2}^kt^{-2r}S_r(y)y^n\\
& & \quad  =-t^{-2k+4}y^nS_{k-1}(y)-t^{-2k+2}x^2y^nS_k(y)+
t^{-2k}y^nS_{k+1}(y)\\
& & \quad -t^2x^2y^n+2t^2x^2\sum_{r=-1}^kt^{-2r}S_r(y)y^n, 
\end{eqnarray*}
and we are done. 
\end{proof}

We now proceed with the proof of Theorem 3.1. From Lemma 3.2 and Lemma
3.3 it follows that 
$A(k,n)=\bar{A}(k,n)$. Also, by rotating the figure by 
$180^\circ$ we see that $A(k,0)=\bar{A}(2p+1-k,0)$. 
In particular, for $1\leq i \leq p$ we have
\begin{eqnarray*}
A(p+i,0)=\bar{A}(p-i+1,0)=A(p-i+1,0).
\end{eqnarray*}
We prove the theorem by induction on $i$. Of course the induction
can start at $i=-1$, since plugging this value into the formula
produces an obvious equality. The case $i=0$ is also easy, so
let  $k\geq 1$ and assume that the formula holds for
$i=-1,0,\ldots, k-1$. From $A(p+k,0)=A(p-k+1,0)$ and Lemma 3.4 it
follows that
\begin{eqnarray*}
& & -t^{-2p-2k+6}S_{p+k-2}(y)-t^{-2p-2k+4}x^2S_{p+k-1}(y)+
t^{-2p-2k+2}S_{p+k}(y)-t^2x^2\\
& & +2t^2x^2\sum_{r=-1}^{p+k-1}t^{-2r}S_r(y)=-t^{-2p+2k+4}S_{p-k-1}(y)
+t^{-2p+2k+2}x^2S_{p-k}{y}\\
& & +t^{-2p+2k}S_{p-k+1}(y)-t^2x^2+2t^2x^2\sum_{r=-1}^{p-k-1}t^{-2r}S_r(y).
\end{eqnarray*}
This yields
\begin{eqnarray*}
& & t^{-2p-2k+2}S_{p+k}(y)+t^{-2p+2k+4}S_{p-k-1}{y}\\
& & =(t^{-2p-2k}+6S_{p+k-2}(y)+t^{-2p+2k}S_{p-k+1}(y))+
x^2(t^{-2p-2k+4}S_{p+k-1}(y)\\
& & +t^{-2p+2k+2}S_{p-k}(y))-2t^2x^2\sum_{r=p-k}^{p+k-1}t^{-2r}S_r(y)\\
& & = t^{-2p+3}[t^{-2k+3}S_{p+k-2}(y)+t^{2k-3}S_{p-k+1}(y)+
x^2(t^{-2k+1}S_{p+k-1}(y)\\
& & +t^{2k-1}S_{p-k}(y))]
-2t^2x^2\sum_{j=0}^{k-1}(t^{-2p-2j}S_{p+j}(y)+t^{-2p+2j+2}S_{p-j-1}(y))\\
& & =t^{-2p+3}[(-1)^{k-2}S_{2k-4}(x)(t^{-1}S_p(y)+tS_{p-1}(y))+
(-1)^{k-1}x^2S_{2k-2}(x)\times \\
& & \times (t^{-1}S_p(y)+tS_{p-1}(y))
 -2x^2\sum _{j=0}^{k-1}(-1)^jS_{2j}(x)(t^{-1}S_p(y)+tS_{p-1}(y))].
\end{eqnarray*}
Now we use the following identity, which can be proved for example
by induction,
\begin{eqnarray*}
2x^2\sum _{r=0}^m(-1)^rS_{2r}(x)=(-1)^{m+1}
(S_{2i-2}(x)-x^2S_{2i}(x)-S_{2i+2}(x)).
\end{eqnarray*}
Applying it we further transform the above expression into
\begin{eqnarray*}
& & t^{-2p+3}(t^{-1}S_p(y)+tS_{p-1}(y))[(-1)^{k-2}S_{2k-4}(x)+
(-1)^{k-1}x^2S_{2k-2}(x)\\
& & -(-1)^k(S_{2k-4}(x)-x^2S_{2k-2}(x)-S_{2k}(x))]\\
& & =(-1)^kt^{-2p+3}S_{2k}(x)(t^{-1}S_p(y)+tS_{p-1}(y)),
\end{eqnarray*}
and the induction is complete. 
All what remains to prove is the case $i=p+1$.
To this end we extend formally the recursive relation
and define $A(0,0)$ and $A(2p+1,0)$. 
From
\begin{eqnarray*}
& & A(2,0)=t^{-2}A(1,1)-t^{-4}A(0,0)+(t^2-t^{-2}+y)x^2\\
 & & A(2p+1,0)=t^{-2}A(2p,1)-t^{-4}A(2p-1,0)+(t^2-t^{-2}+y)x^2
\end{eqnarray*}
and  the formulas  for $A(2,0),A(1,1), A(2p,1) $ and
$A(2p-1,0)$ we obtain
\begin{eqnarray*}
A(2p+1,0)=A(0,0)=t^6+t^2-t^2x^2.
\end{eqnarray*}
But this is all that we need for  the above induction argument to work
up to the case $i=p+1$, and the theorem is proved.

\section{Two Curves on the Boundary}

For proving Theorem 1.1 we need to determine the images 
in $K_t(M_p)$ through the inclusion map 
 of the skeins $(1,k)_T$ from the skein module of
the boundary torus. The particular case
of the product-to-sum formula
\begin{eqnarray*}
(0,1)_T*(1,k)_T=t^{-1}(1,k+1)_T+t(1,k-1)_T
\end{eqnarray*}
becomes after projecting to $K_t(M_p)$ 
\begin{eqnarray*}
x\pi((1,k)_T)=t^{-1}\pi((1,k+1)_T)+t\pi((1,k-1)_T),
\end{eqnarray*}
where again, the multiplication on the left-hand side means that
after we expand $\pi((1,k)_T)$ in terms of basis elements, the
powers of $x$ are raised by $1$. 
Hence we have a recursive formula for computing $(1,k)_T$. 

This shows that it suffices to compute $(1,k)_T$ and $(1,k-1)_T$ for some 
value of $k$. 
In the following two propositions the reasoning will be done on figures
depicting  the $(2,5)$-torus knot.

{\bf Proposition~4.1.} {\it
For all $p\in {\mathbb Z}$,
\begin{eqnarray*}
\pi((1,-4p-2)_T)=2+(-1)^{p+1}\left(t^{-2p-3}S_{2p+2}(x)-
t^{-2p+1}S_{2p-2}(x)\right)\\ \times (tS_{p-1}(y)+t^{-1}S_p(y)).
\end{eqnarray*}
}

\begin{proof}
The curve $(1,-4p-2)_T$ is depicted with dotted line  in Fig. 4.1.

\begin{figure}[htbp]
\centering
\leavevmode
\epsfxsize=2.2in
\epsfysize=1.2in
\epsfbox{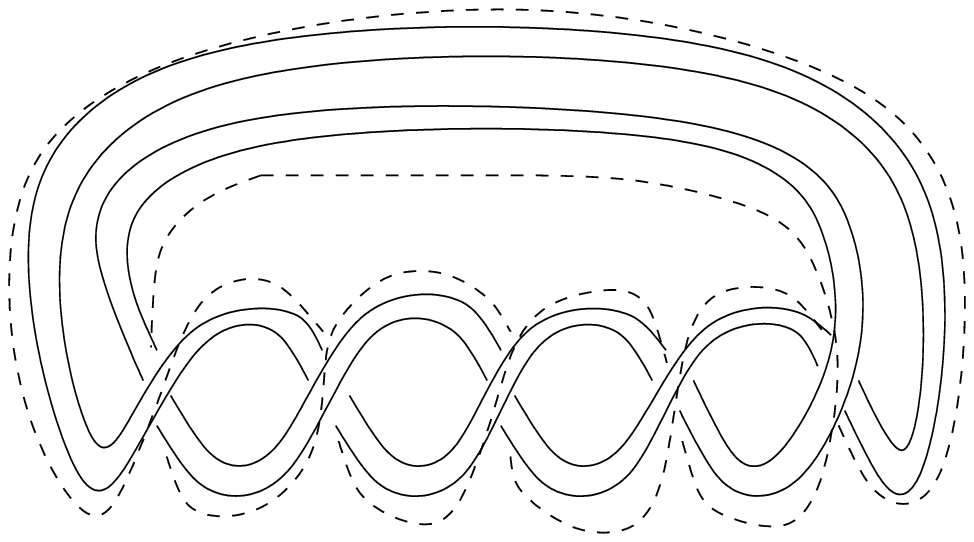}

Figure 4.1.    
\end{figure}

The computation of
$\pi((1,-4p-2)_T)$ is simplified by the following observation.
If we remove from $M_p$ a regular neighborhood of the M\"{o}bius band
that the knot bounds, then the resulting 3-manifold N
 is a solid torus that contains
$\pi((1,-4p-2)_T)$ in its interior. So if we express the skein 
$\pi((1,-4p-2)_T)$ in terms of the basis of $K_t(M_p)$ only powers
of $y$ should appear.
The coefficients of these powers are
computed as for the solid torus $N$ and the image
of the curve $(2p+1,2)_T$ in its boundary. Note that when we turn the
torus inside out the meridian and the longitude are exchanged.
  Using the formula deduced in Section 5 of \cite{FG}, we find that
this image, whether considered in $K_t(M_p)$ or $K_t(N)$, is 
\begin{eqnarray*}
(-1)^{2}t^{-2(2p+1)}\left(t^{-4}S_{2p+1}(y)-t^{4}S_{2p-1}(y)\right).
\end{eqnarray*}
Using Theorem 3.1 for $i=p+1$ and $i=p-1$
 we deduce that this is further equal to
\begin{eqnarray*}
 & & -S_{-2}(y) +(-1)^{p+1}S_{2p+2}(x)(tS_{p-1}(y)+t^{-1}S_p(y))\\
&  & \quad +
S_0(y)-(-1)^{p-1}t^{-2p+1}S_{2p-2}(x)(tS_{p-1}(y)+t^{-1}S_p(y))\\
& & =2+(-1)^{p+1}(t^{-2p-3}S_{2p+2}(x)-
t^{-2p+1}S_{2p-2}(x))(tS_{p-1}(y)+t^{-1}S_p(y)).
\end{eqnarray*}
and we are done. 
\end{proof}

For expanding the second boundary curve in terms of the basis we will
need the following formula.

{\bf Lemma~4.2.} {\it
For any integer $p$ the following holds
\begin{eqnarray*}
\sum_{k=1}^{p-1}S_{2p-4i-2}=-\frac{S_{2p-3}}{S_1}.
\end{eqnarray*}
}

\begin{proof}
First, since $S_{n}=-S_{-n-2}$, 
\begin{eqnarray*}
\sum_{k=1}^{p-1}S_{2p-4i-2}=\sum_{k=1}^{p-1}(-1)^kS_{2p-2k-2}.
\end{eqnarray*}
Recall that $S_n(2\cos x)=\sin (n+1)x/\sin x$ and that
\begin{eqnarray*}
\sum_{k=1}^n \sin (2k+1)x= \sin(2n+2)x/2\cos x.
\end{eqnarray*}  
We have
\begin{eqnarray*}
& & 
\sum_{k=1}^{p-1}(-1)^kS_{2p-2k-2}(2\cos x)=\sum_{k=1}^{p-1}(-1)^k\frac{\sin(
2p-2k-1)x}{\sin x}\\
& & \quad =-\frac{\sin (2p-2)x}{2\sin x\cos x}=-\frac{\sin (2p-2)x}{\sin 2x}\\
& & \quad =
-\frac{\sin (2p-2)x}{\sin x}\cdot \frac{\sin x}{\sin 2x}=- \frac{S_{2p-3}(2\cos
x)}{S_1(2\cos x)}
\end{eqnarray*}
and the identity follows.
\end{proof}

{\bf Proposition~4.3.} {\it
For all $p\in {\mathbb Z}$,
\begin{eqnarray*}
\pi((1,-4p-1)_T)=tm+(-1)^{p+1}\left(t^{-2p-2}S_{2p+1}(x)-t^{-2p+2}S_{2p-3}(x)
\right) \\ \times (tS_{p-1}(y)+
t^{-1}S_p(y)).
\end{eqnarray*}
}

\begin{proof}
To make our figures easier to read we replace the drawing of 
the $(2,2p+1)$-torus knot as described in Fig. 4.2, and ask the reader
to use the imagination and think of the 3-valent graph as being the 
torus knot. 

\begin{figure}[htbp]
\centering
\leavevmode
\epsfxsize=3.3in
\epsfysize=1.1in
\epsfbox{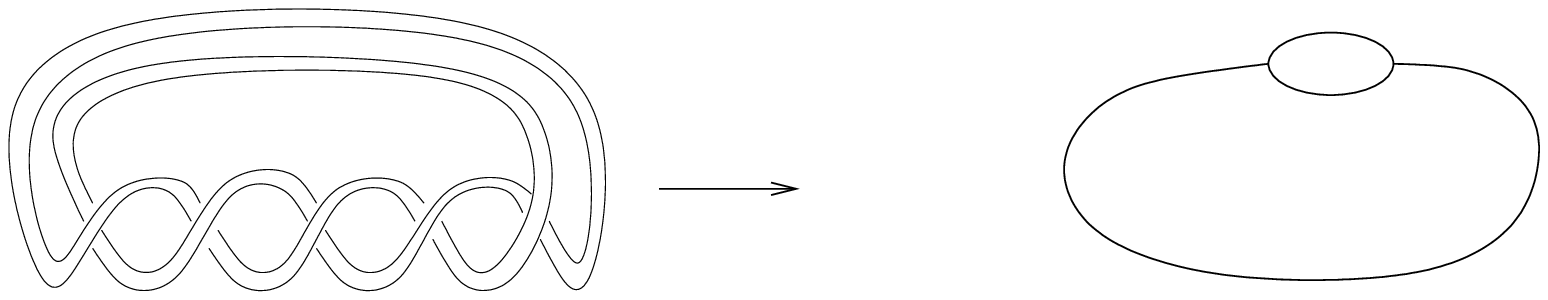}

Figure 4.2.    
\end{figure}

As such,  the dotted curve from Fig. 4.3  multiplied by
$t^{-6p}$ represents the skein  $\pi((1,-4p-1)_T)$. The reason for the 
multiplcation by  a power of $t$ is that, as drawn, the curve has
the blackboard framing, while the framing of $(1,-4p-1)_T$ isparallel
to the boundary torus. Let us write the skein from Fig. 4.3 
in terms of the basis of the module  $K_t(M_p)$. 

\begin{figure}[htbp]
\centering
\leavevmode
\epsfxsize=2.2in
\epsfysize=1.2in
\epsfbox{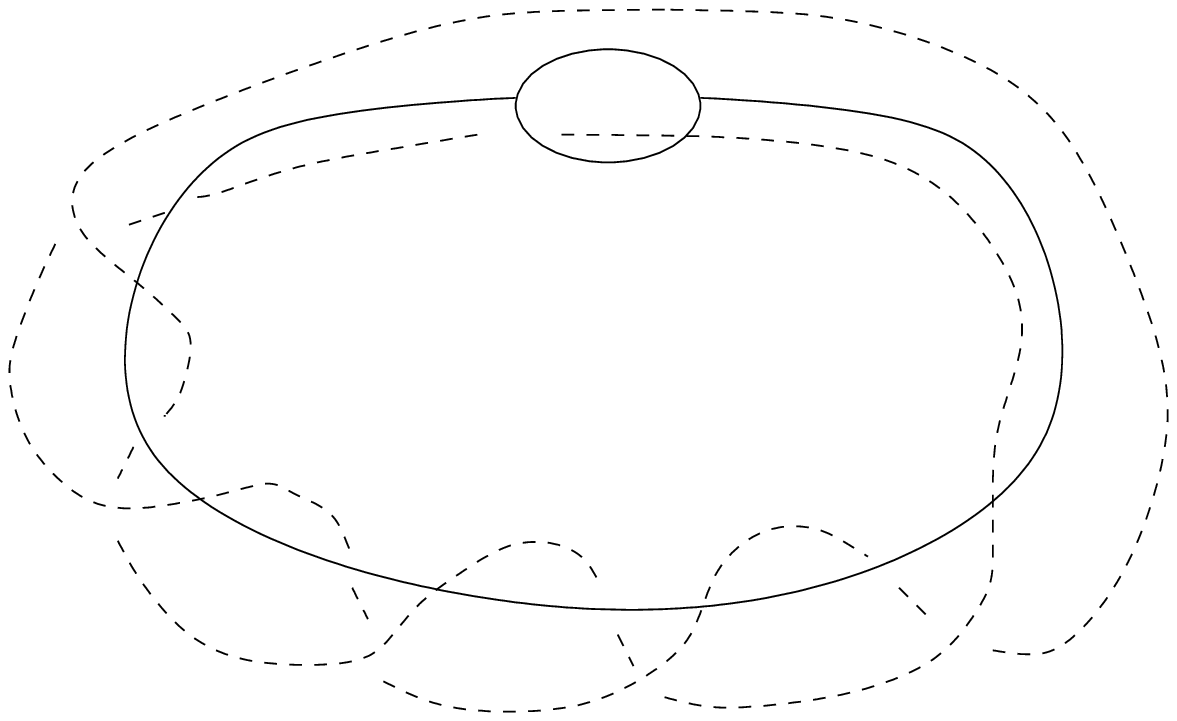}

Figure 4.3.    
\end{figure}

We consider the skeins $a_k$ and $b_k$ described in
Fig. 4.4  where $a_k$ has  $k$ self crossings, $k\leq 2p+1$. 
Resolving the rightmost crossing of $a_k$ and performing the 
twistings we obtain the recursive
relations
\begin{eqnarray*}
&  & a_{2k}=t(-t^3)^{2k-1}b_{2p-2k+1}+t^{-1}a_{2k-1}\\
&  & a_{2k+1}=t(-t^3)^{2k}b_{2p-2k+1}+t^{-1}a_{2k}.
\end{eqnarray*}
Also, $a_0=(-t^2-t^{-2})b_{2p+1}$. Let us determine $b_{k}$ in terms of
$x$ and $y$. 

\begin{figure}[htbp]
\centering
\leavevmode
\epsfxsize=5in
\epsfysize=1.4in
\epsfbox{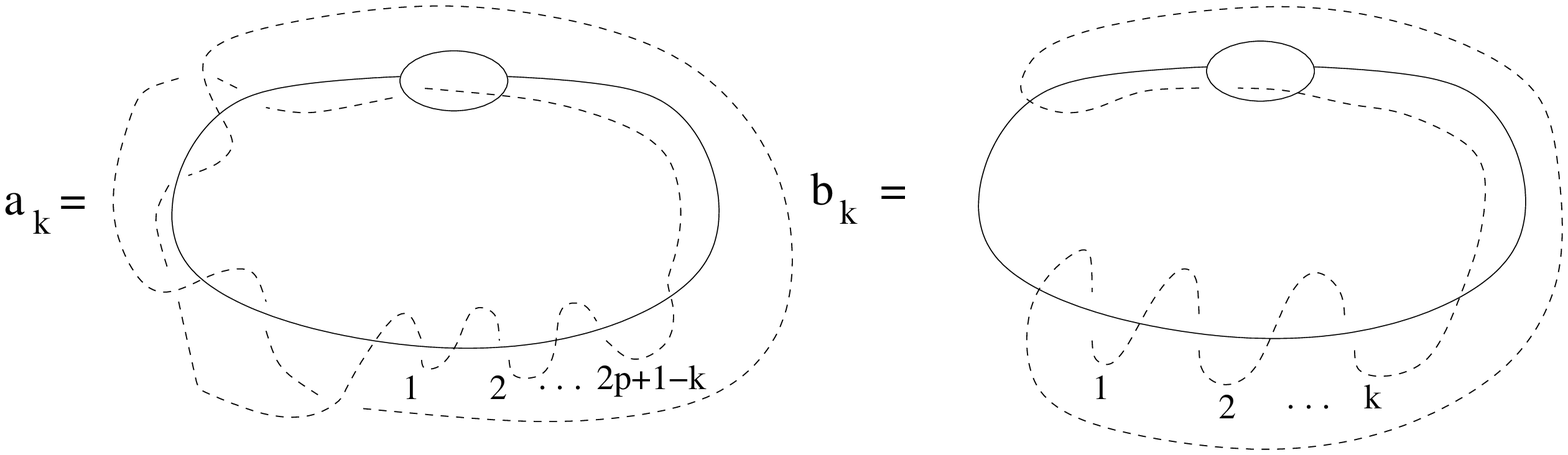}

Figure 4.4. 
   
\centering
\leavevmode
\epsfxsize=4.5in
\epsfysize=1.2in
\epsfbox{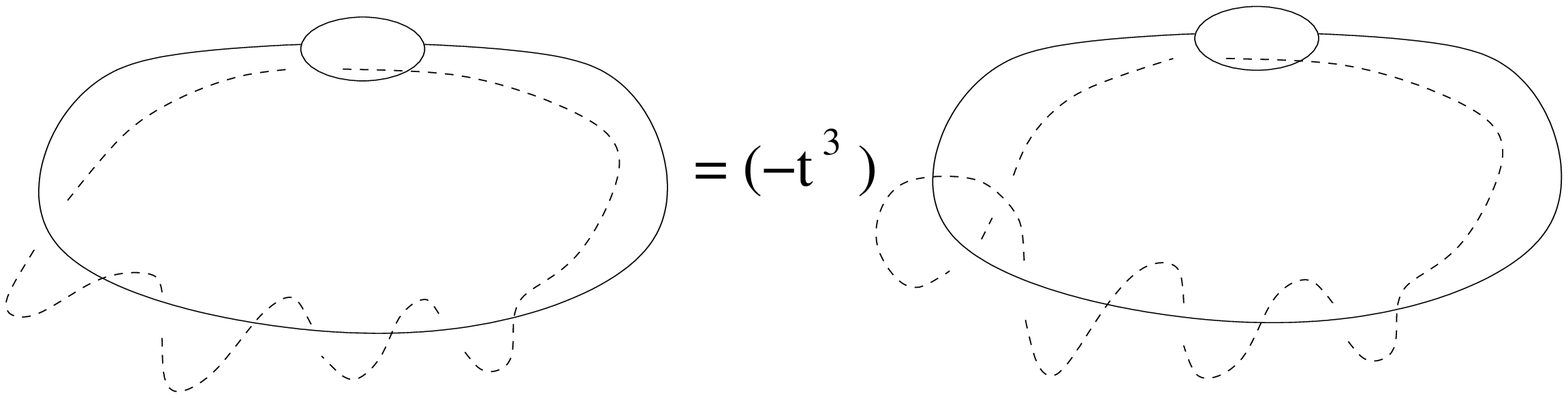}

Figure 4.5. 
\end{figure}

Transform $b_{k}$ as in Fig. 4.5, then resolve the crossing to 
obtain the recursive relation
\begin{eqnarray*}
b_{k+1}=(-t)^3(tb_{k-1}+t^{-1}b_{k-2}y)=-t^2yb_k-t^4b_{k-1}.
\end{eqnarray*}
If we let $x_k=(-t^{-2})^kb_k$, then $x_{k+1}=yx_k-x_{k-1}$, and
$x_0=m$, $x_1=t^{-2}m$. Solving the second order recurrence we
obtain $x_k=-t^{-2}xS_{k-1}(y)-xS_{j-2}(y)$, hence 
\begin{eqnarray*} 
b_k=(-1)^{k+1}t^{2k-2}xS_{k-1}(y)+(-1)^{k+1}t^{2k}xS_{k-2}(y), \quad \mbox{ for
all } k\geq 0.
\end{eqnarray*}
Let us compute $a_{2p+1}$ by writing the solution for the  nonhomogeneous
first order  recursive relation that we deduced for $a_k$. We have
\begin{eqnarray*}
&  & a_{2p+1}=t^{-2p-1}a_0+t\sum _0^p t^{-2i}(-t^3)^{2p-2i}b_{2i+1}+
t\sum _{1}^p t^{-2i+1}(-t^3)^{2p-2i+1}b_{2i-1}\\
& & \quad  = t^{-2p-1}a_0+t\sum _0^p t^{-2i}(-t^3)^{2p-2i}b_{2i+1}+
t\sum _{0}^{p-1} t^{-2i-1}(-t^3)^{2p-2i-1}b_{2i+1}\\
& & \quad 
=  t^{-2p-1}(-t^2-t^{-2})b_{2p+1}+t^{-2p+1}b_{2p+1}+t^{6p+1}\sum_0^{p-1}
t^{-8i}(1-t^{-4})b_{2i+1}\\
& & \quad = -t^{2p-3}b_{2p+1}+t^{6p+1}\sum_0^{p-1}t^{-8i}(1-t^4)b_{2i+1}.
\end{eqnarray*}
If we replace in this the value we computed for $b_{k}$ we obtain
\begin{eqnarray*}
& & a_{2p+1}=-t^{-2p-3}\left[t^{4p}xS_{2p}(y)+t^{4p+2}xS_{2p-1}(y)\right]\\
& & +t^{6p+1}(1-t^{-4})x\sum_{0}^{p-1}\left[t^{-4i}S_{2i}(y)+t^{-4i+2}S_{2i-1}
(y)\right].
\end{eqnarray*}
Factoring out an $t^{1-2p}$ and changing 
the order of summation for the second of the two sums we obtain 
that this is further equal to
\begin{eqnarray*}
& & -t^{-2p-3}\left(t^{4p}xS_{2p}(y)+t^{4p+2}xS_{2p-1}(y)\right)\\
& & \quad 
+t^{4p+2}(1-t^{-4})x\sum_{0}^{p-1}\left[t^{2p-4i+1}S_{2i}(y)+t^{-2p+4i+1}
S_{p+(p-1-2i)}(y)\right].
\end{eqnarray*}
Note that $p-(p-1-2i)-1=2i$ hence we can replace each square bracket using 
Theorem 3.1. We also replace $S_{2p}(y)$ and $S_{2p-1}(y)$ to
 deduce that
\begin{eqnarray*}
& & a_{2p+1}=-t^{2p-3}x\left[(-1)^pt^{2p+1}S_{2p}(x)(tS_{p-1}(y)+t^{-1}S_p(y))
\right]\\
& & \quad -t^{2p-1}x\left[-t^{4p-2}+(-1)^{p-1}t^{2p-1}S_{2p-2}(x)(tS_{p-1}(y)+
t^{-1}S_{p}(y))\right]\\
& & \quad + t^{4p+2}(1-t^{-4})x\sum _{i=0}^{p-1}(-1)^{p-1-2i}S_{2(p-1-2i)}(x)
(tS_{p-1}(y)+t^{-1}S_p(y)).
\end{eqnarray*}
Since $S_{2p}-S_{2p-2}=T_{2p}$ we conclude that
\begin{eqnarray*}
a_{2p+1}=t^{6p+1}x+(-1)^{p+1}\left[t^{4p-2}xT_{2p}(x)+
(t^{4p+2}-t^{4p-2})x\sum_{0}^{p-1}S_{2p-4i-2}(x)\right]\\
\times (tS_{p-1}(y)+t^{-1}
S_p(y)).
\end{eqnarray*}
But
 $\pi((1,-4p-1)_T)=t^{-6p}a_{2p+1}$, so
\begin{eqnarray*}
& & \pi((1,-4p-1)_T)=tx+(-1)^{p+1}[t^{-2p-2}xT_{2p}(x)\\ & & 
\quad +
(t^{-2p+2}-t^{-2p-2})x\sum _1^{p-1}S_{2p-4i-2}(x)](tS_{p-1}(y)+
t^{-1}S_p(y)).
\end{eqnarray*}
Applying Lemma 4.2 we get
\begin{eqnarray*}
& &  \pi((1,-4p-1)_T)=tx+(-1)^{p+1}[t^{-2p-2}xT_{2p}(x)\\ & & \quad 
 -
(t^{-2p+2}-t^{-2p-2})S_{2p-3}(x)](tS_{p-1}(y)+
t^{-1}S_p(y)),
\end{eqnarray*}
and the desired formula follows from
$xT_{2p}(x)=T_{2p+1}(x)+T_{2p-1}(x)=S_{2p+1}(x)-S_{2p-3}(x)$.
\end{proof} 

Combining the two propositions we obtain

{\bf Proposition~4.4.} {\it
For all $k\in{\mathbb Z}$ one has
\begin{eqnarray*}
& & 
\pi((1,k)_T)=t^{4p+k+2}T_{4p+k+2}(x)+(-1)^{p+1}(t^{2p+k-1}S_{-2p-k}(x)\\ & &
\quad -
t^{2p+k+3}S_{-2p-k-4}(x))\times (tS_{p-1}(y)+t^{-1}S_p(y)).
\end{eqnarray*}
}

\begin{proof}
The particular case of the  product-to-sum formula
\begin{eqnarray*}
(1,k+1)_T=t(0,1)_T*(1,k)_T-t^2(1,k-1)_T
\end{eqnarray*}
descends in the knot complement to the recursive relation
\begin{eqnarray*}
\pi((1,k+1)_T)=tx\pi((1,k)_T)-t^2\pi((1,k-1)_T).
\end{eqnarray*}
Hence the values of all $\pi((1,k)_T$ can be computed from
$\pi((1,-4p-2)_T)$ and $\pi((1,-4p-1)_T)$.  An easy induction
proves the proposition. 
\end{proof}

\section{Proof of the Main Result}

We start by constructing an element in the peripheral ideal of the 
$(2,2p+1)$-torus knot using the formulas from Section 4.  

{\bf Lemma~5.1.} {\it
The element 
\begin{eqnarray*}
(1,-2p-3)_T-t^{-8}(1,-2p+1)_T+t^{2p-5}(0,2p+3)_T-
t^{2p-1}(0,2p-1)_T
\end{eqnarray*}
is in the peripheral ideal of the $(2,2p+1)$-torus knot.
}

\begin{proof}
Applying Proposition 4.4 and using the fact that $S_{-1}=0$  we obtain 
\begin{eqnarray*}
\pi((1,-2p-3)_T)=t^{2p-1}T_{2p-1}(x)+(-1)^{p+1}(t^{-4}S_3(x)(tS_{p-1}(y)+
t^{-1}S_p(y))
\end{eqnarray*}
and 
\begin{eqnarray*}
& & 
\pi((1,-2p+1)_T)=t^{2p+3}T_{2p+3}(x)+(-1)^{p+1}(-t^{4}S_{-5}(x))(tS_{p-1}(y)+
t^{-1}S_p(y))\\
& & =t^{2p+3}T_{2p+3}(x)+(-1)^{p+1}t^{4}S_{3}(x)(tS_{p-1}(y)+
t^{-1}S_p(y)).
\end{eqnarray*}
Multiplying the second equality by $t^{-8}$ and subtracting it from the
first we obtain
\begin{eqnarray*}
(1,-2p-3)_T-t^{-8}(1,-2p+1)_T=t^{2p-1}T_{2p-1}(x)-t^{2p-5}T_{2p+3}(x).
\end{eqnarray*}
But $T_n(x)=\pi((0,n)_T)$ for all $n$. Hence the image through $\pi $ of
\begin{eqnarray*}
(1,-2p-3)_T-t^{-8}(1,-2p+1)_T+t^{2p-5}(0,2p+3)_T-
t^{2p-1}(0,2p-1)_T
\end{eqnarray*}
is zero. This shows that this element is in the kernel of $\pi$, that is
in the peripheral ideal. 
\end{proof}

{\bf Proposition~5.2.} {\it
The polynomial 
\begin{eqnarray*}
(l-t^{-4}lm^4+t^{4p-2}-t^{4p+10}m^4)(l-t^{4p+2}m^{4p+2})
\end{eqnarray*}
is in the noncommutative A-ideal of the $(2,2p+1)$-torus knot.\\
}

\begin{proof}
Through the inclusion 
\begin{eqnarray*}
K_t({\mathbb T}^2\times I)\hookrightarrow {\mathbb C}_t[l,l^{-1},m,m^{-1}]
\end{eqnarray*}
an element $(a,b)_T$ transforms into $t^{-ab}(l^am^b+l^{-a}m^{-b})$. 
Applying this to the element provided by Lemma 5.1 we deduce that
the extension of the peripheral ideal $I_t(K)$
 to ${\mathbb C}_t[l,l^{-1},m,m^{-1}]$
contains the element
\begin{eqnarray*}
& & t^{2p+3}lm^{-2p-3}+t^{2p+3}l^{-1}m^{2p+3}-t^{2p-9}lm^{-2p+1}-t^{2p-9}l^{-1}
m^{2p-1}\\
& & +t^{2p-5}m^{2p+3}+t^{2p-5}m^{-2p-3}-t^{2p-1}m^{2p-1}-t^{2p-1}m^{-2p+1}.
\end{eqnarray*}
Now contract the extension $I_t(K)$ to ${\mathbb C}_t[l,m]$. The above
element gives rise through multiplication to the left 
 by $lm^{2p+3} $  
to a polynomial in the noncommutative A-ideal. This polynomial
is
\begin{eqnarray*}
& & t^{-2p-3}l^2+t^{6p+9}m^{4p+6}-t^{-2p-14}l^2m^4
 -t^{6p-3}m^{4p+2}+t^{2p-5}lm^{4p+6}\\ & & +t^{2p-5}l
 -t^{2p-1}lm^{4p+2}-t^{2p-1}lm^4.
\end{eqnarray*}
After multiplying by $t^{2p+3} $ and factoring this becomes the polynomial
from the statement.
\end{proof}

Let us proceed with the proof of Theorem 1.1. To this end we rephrase
slightly Theorem 2 in \cite{G2}:

\medskip
\begin{em}
Let $K$ be a knot whose A-ideal  ${\mathcal A}_t(K)$
contains a polynomial $\sum _{p,q}\gamma _{p,q}l^pm^q$
 of degree $2$ in $l$ such that there exists no $n\geq 0$ for which the 
expression 
$\sum_{q}\gamma _{2,q}(-1)^qt^{(2n+2)q}$ is identically equal to zero.
Assume in addition that this polynomial arises 
from a skein in the peripheral ideal of the knot. Then for any knot 
$K'$ with the property that ${\mathcal A}_t(K)={\mathcal A}_t(K')$,
one has  $\kappa_n(K)=\kappa_n(K')$ for all $n\geq 1$.
\end{em}
\medskip

Here the assumption that the polynomial comes from a skein in the
peripheral ideal yields to a simpler version of the relation among
the coefficients of the polynomial to be checked
than the one given in \cite{G2},
and the same proof works. Recall also that we denoted
by $\kappa _n(K)$ the $n$th colored Kauffman bracket of $K$, which
is (up to a change of variable) the same as the $n$th colored Jones polynomial
as defined in \cite{RT}.

Expand the polynomial  from Proposition 5.2 to get
\begin{eqnarray*}
& & l^2+t^{8p+12}m^{4p+6}-t^{-12}l^2m^4 -t^{8p}m^{4p+2}+t^{4p-2}lm^{4p+6}\\
& & +t^{4p-2}l -t^{4p+2}lm^{4p+2}-t^{4p+2}lm^4.
\end{eqnarray*}
The only coefficients $\gamma _{2,q} $ that are not zero are
$\gamma _{2,0}=1$ and $\gamma _{2,4}=-t^{-12}$. The expression
\begin{eqnarray*}
1-t^{-12}t^{4(2n+2)}
\end{eqnarray*} 
is identically  zero if and only if $8n-4=0$ for some $n$.
This has to be  so since  $t$ is the variable of a polynomial.
The equality cannot hold, which  proves Theorem 1.1.  

Through the procedure described in \cite{G2} we can find the following 
recursive relation for the colored Kauffman brackets of the $(2,2p+1)$-torus
knot
\begin{eqnarray*}
& & (-t^{-4np-6n-6p-9}+t^{-4np+2n-6p-5})\kappa_{n+1}(K)\\ & & +
(-t^{4np+6n+6p+1}-t^{-4np-6n-2p-11}+t^{4np-2n+6p-3}+
t^{-4np+2n-2p+1})\kappa _n(K)\\
& & +(t^{-4np-6n+2p+3}-
t^{4np+6n+6p+9}-t^{-4np+2p+2n-9}+t^{4np-2n+2p-11})\kappa_{n-1}(K)\\
& & +(t^{4np+6n-2p-11}+t^{-4np-6n+6p+1}-
t^{4np-2n-2p+1}-t^{-4np+2n+6p-3})\kappa _{n-2}(K)\\
& & +(t^{4np+6n-2p-3}-t^{4np-2n-6p-7})\kappa_{n-3}(K)=0.
\end{eqnarray*}
Recall that $\kappa _{-3}(K)=-\kappa _1(K)$, $\kappa _{-2}(K)=-1$, 
$\kappa _{-1}(K)=0$, $\kappa _0(K)=1$, so this recurrence 
 allows us indeed to compute all colored Kauffman brackets.

\end{document}